\documentclass[a4paper, 11pt]{article}
\usepackage{amsmath, amsfonts, amssymb}
\usepackage{amscd}
\usepackage[english]{babel}
\newtheorem{Th}{Theorem}

\newtheorem{Prop}{Proposition}

\newtheorem{Ex}{Example}
\newenvironment{Proof} 
{\par\noindent{\bf Proof.}} 
{\hfill$\scriptstyle\blacksquare$}

\begin {document}

\begin{center}
\textbf {On metrisation of the space of idempotent probability measures} \\

\medskip
\textbf {Adilbek Atakhanovich Zaitov}\\
\smallskip
{Tashkent institute of architecture and civil engineering}\\
https://orcid.org/0000-0002-2248-0442
\end{center}

\begin{abstract}
In this paper we introduce a metrics on the space of idempotent probability measures on a given compactum, which extends the metrics on the compactum. It is proven the introduced metrics generates the pointwise convergence topology on the space of idempotent probability measures.\\


\end{abstract}

\begin{flushright}
\textit{Dedicated to the kind memory of\\
professor Fedorchuk, Vitalii Vital'evich\\
with admiration and gratitude}
\end{flushright}

\section{Itroduction}

Idempotent mathematics is based on replacing the usual arithmetic operations with a new set of basic operations, i.~e., on replacing numerical fields by idempotent semirings and semifields. Typical example is given by the so-called $\max$-$\mbox{plus}$ algebra $\mathbb{R}_{\max}$. Let $\mathbb{R}$ be the field of real numbers. Then $\mathbb{R}_{\max} = \mathbb{R} \cup \{-\infty\}$ with the operations $x \oplus y = \max\{x,\, y\}$ and $x \odot y = x + y$. The new addition $\oplus$ is idempotent, i.~e. $x \oplus x = x$ for all elements $x$.

Many authors (S.~C.~Kleene, S.~N.~N.~Pandit, N.~N.~Vorobjev, B.~A.~Carr\'{e}, R.~A.~Cuninghame-Green, K.~Zimmermann, U.~Zimmermann, M.~Gondran, F.~L.~Baccelli, G.~Cohen, S.~Gaubert, G.~J.~Olsder, J.-P.~Quadrat, and others) used idempotent semirings and matrices over these semirings for solving some applied problems in computer science and discrete mathematics, starting from the classical paper by S.~C.~Kleene \cite{Kleene1956}.

The modern idempotent analysis (or idempotent calculus, or idempotent mathematics) was founded by V.~P.~Maslov and his
collaborators \cite{Litv2007}.  Some preliminary results are due to E.~Hopf and G.~Choquet, see \cite{Choquet1955}, \cite{Hopf1950}.

Idempotent mathematics can be treated as the result of a dequantization of the traditional mathematics over numerical fields as the Planck constant $h$ tends to zero taking imaginary values. This point of view was presented by G.~L.~Litvinov and V.~P.~Maslov \cite{LMSh2002}. In other words, idempotent mathematics is an asymptotic version of the traditional mathematics over the fields of real and complex numbers.

The basic paradigm is expressed in terms of an idempotent correspondence principle. This principle is closely related to the well-known correspondence principle of N.~Bohr in quantum theory. Actually, there exists a heuristic correspondence between important, interesting, and useful constructions and results of the traditional mathematics over fields and analogous constructions and results over idempotent semirings and semifields (i.~e., semirings and semifields with idempotent addition).

A systematic and consistent application of the idempotent correspondence principle leads to a variety of results, often quite unexpected. As a result, in parallel with the traditional mathematics over fields, its ``shadow,'' idempotent mathematics, appears. This ``shadow'' stands approximately in the same relation to traditional mathematics as classical physics does to quantum theory.

Remind \cite{Litv2007} a set $S$ equipped with two algebraic operations: addition $\oplus$ and multiplication $\odot$, is said to be a \textit{semiring} if the following conditions are satisfied:
\begin{itemize}
\item the addition~$\oplus$ and the multiplication~$\odot$ are associative;
\item the addition~$\oplus$ is commutative;
\item the multiplication~$\odot$ is distributive with respect to the addition~$\oplus$:
\begin{gather*}
x\odot(y\oplus z) = x\odot y \oplus x\odot z \quad\mbox{ and }\\
(x\oplus y)\odot z = x\odot z \oplus y\odot z
\end{gather*}
for all ~$x,\ y,\ z\in S$.
\end{itemize}

A \textit{unit} of a semiring~$S$ is an element~$\textbf{1}\in S$ such that~$\textbf{1}\odot x = x\odot \textbf{1} = x$  for all~$x\in S$.
A \textit{zero} of a semiring~$S$ is an element~$\textbf{0}\in S$ such that~$\textbf{0}\ne \textbf{1}$ and ~$\textbf{0}\oplus x = x\oplus \textbf{0} = x$  for all~$x\in S$. A semiring~$S$ is called an \textit{idempotent semiring} if~$ x\oplus x = x$ for all~$x\in S$. A (an idempotent) semiring~$S$ with neutral elements~$\textbf{0}$ and~$\textbf{1}$ is called a (an \textit{idempotent}) \textit{semifield} if every nonzero element of~$S$ is invertible. Note that dio\"{\i}ds, quantales and inclines are examples of idempotent semirings \cite{Litv2007}.

Let us state \textit{Maslov dequantization}. Let~$\mathbb{R}=(-\infty,\, +\infty)$ be the field of real numbers and~$\mathbb{R}_+=[0,\, +\infty)$ be the semiring of all nonnegative real numbers (with respect to the usual addition and multiplication). Consider a map~$\Phi_h\colon \mathbb{R}_{+}\to S=\mathbb{R}\cup\{-\infty\}$ defined by the equality
\begin{gather*}
\Phi_h(x)= h\ln\, x, \qquad h>0.
\label{1.1}
\end{gather*}
Let us undergo the usual operations of addition and multiplication from~$\mathbb{R}_{+}$ into~$S$ using the map~$\Phi_h$. Let
\begin{gather*}
u=\Phi_h(x)= h\ln\, x,\qquad v=\Phi_h(y)= h\ln\, y.
\end{gather*}
Then
\begin{gather*}
\Phi_h(x+y)=h\ln\,(x+y) = h\ln\,\left(\textrm{e}^{\frac{u}{h}}+\textrm{e}^{\frac{v}{h}}\right),\\
\Phi_h(xy)=h\ln\,(xy)=h\ln\, x+h\ln\, y.
\end{gather*}
Put~$u\oplus_h v = \Phi_h(x+y)$ and~$u\odot v = \Phi_h(xy)$, i.~e.~$u\oplus_h v = h\ln\, \left(\textrm{e}^{\frac{u}{h}}+\textrm{e}^{\frac{v}{h}}\right)$ and~$u\odot v = u+v$. The imagine~$\Phi_h(0)=-\infty$ of the usual zero~$0$ is a zero~$\textbf{0}$ and the imagine~$\Phi_h(1)=0$ of the usual unit~$1$ is a unit~$\textbf{1}$ in~$S$ with respect to these new operations. Thus~$S$ obtains the structure of a semiring~$\mathbb{R}^{(h)}$ isomorphic to~$\mathbb{R}_+$.

The direct check shows that~$u\oplus_h v\to \max\{u,\, v\}$ as~$h\to 0$. The convention $-\infty\odot x =-\infty$ allows us to extend $\oplus$ and $\odot$ over $S$. It can easily be checked that~$S$ forms a semiring with respect to the addition~$u\oplus v=\max\{u,\, v\}$ and the multiplication~$u\odot v=u+v$ with zero~$\textbf{0}=-\infty$ and unit~$\textbf{1}=0$. Denote this semiring by~$\mathbb{R}_{\max}$; it is idempotent, i.~e., $u \oplus u =u$ for all its elements. The semiring~$\mathbb{R}_{\max}$ is actually a semifield. The analogy with quantization is obvious; the parameter $h$ plays the role of the Planck constant, so $\mathbb{R}_+$ can be
viewed as a ``quantum object'' and $\mathbb{R}_{\max}$ as the result of its ``dequantization''. This passage to $\mathbb{R}_{\max}$ is called the \textit{Maslov dequantization}.

The notion of idempotent (Maslov) measure finds important applications in different part of mathematics, mathematical
physics and economics (see the survey article \cite{Litv2007} and the bibliography therein). Topological and categorical
properties of the functor of idempotent measures were studied in \cite{ZaitKhol2014}, \cite{Zar2010}. Although idempotent measures are not additive and corresponding functionals are not linear, there are some parallels between topological properties of the
functor of probability measures and the functor of idempotent measures (see, for example \cite{ZaitKhol2014}) which are
based on existence of natural equiconnectedness structure on both functors.

However, some differences appear when the problem of the metrisability of the space of idempotent probability measures was studying. The
problem of the metrisability of the space of the usual probability measures was investigated in \cite{Fedorchuk1990}. We show that the analog of the metrics introduced in \cite{Fedorchuk1990} (on the space of probability measures) is not metrics on  the space of idempotent probability measures. We show the mentioned analog is only a pseudometrics.

In this paper we  introduce a metrics on the space of idempotent probability measures.\\

\section{Idempotent probability measures. Preliminaries}

Let $X$ be a compact Hausdorff space ($\equiv$ a compact), $C(X)$ be the algebra of continuous functions on $X$ with usual algebraic
operations. On $C(X)$ operations $\oplus $  and $\odot $ we will determine by rules $\varphi \oplus \psi =\max \{\varphi ,\psi \}$ and $\varphi \odot \psi=\varphi + \psi$ where $\varphi $, $\psi \in C(X)$.

Remind a functional $\mu :C(X)\to \mathbb{R}$ is called $\left[ 7 \right]$ to be an idempotent probability measure on $X$ if
it satisfies the following properties:

(1) $\mu ({{\lambda }_{X}})=\lambda $ for all $\lambda \in \mathbb{R}$, where  ${{\lambda }_{X}}$ -- constant function;

(2) $\mu(\lambda \odot \varphi)=\lambda \odot \mu(\varphi)$ for all $\lambda \in \mathbb{R}$ и $\varphi \in C(X)$;

(3) $\mu (\varphi \oplus \psi )=\mu (\varphi )\oplus \mu (\psi )$ for all $\varphi $, $\psi \in C(X)$.

For a compact $X$ we denote by $I(X)$ the set of all idempotent probability measures on $X$.
$I(X)$ is a subset of ${{\mathbb{R}}^{C(X)}}$. Really, since $\varphi\oplus\psi=\psi$ for any pair $\varphi, \psi\in C(X)$ with $\varphi\le\psi$ we have $\mu(\varphi)\le \mu(\varphi)\oplus\mu(\psi)= \mu(\varphi\oplus\psi)= \mu(\psi)$, i.~e. $\mu$ is order-preserving functional. That is why $\mu \in \prod\limits_{\varphi\in C(X)}\left[-\|\varphi\|,\ \|\varphi\|\right]$. We consider $I(X)$ as a subspace of ${{\mathbb{R}}^{C(X)}}$.
Sets of the view
\begin{gather*}
\left\langle\mu;\ \varphi_1,\, \dots,\, \varphi_n;\ \varepsilon\right\rangle = \{\nu\in I(X):\ |\nu(\varphi_i)-\mu(\varphi_i)|<\varepsilon,\ i=1,\, \dots,\, n\}
\end{gather*}
where $\varphi_i\in C(X)$, $i=1,\, \dots,\ n$, and $\varepsilon>0$, form a base of open neighbourhoods of given idempotent probability measure $\mu\in I(X)$ according to induced topology.

Let $X$, $Y$ be compacts and $f:X\to Y$ be a continuous map. It is easy to check that a map $I(f):I(X)\to I(Y)$ determined by the formula $I(f)(\mu )(\psi )=\mu (\psi \circ f)$ is continuous. The construction  $I$ is a normal functor acting in the category compacts and their continuous maps. Therefore for each idempotent probability measure  $\mu \in I(X)$ one may determine its \textit{support}:

\begin{gather*}
\text{supp}\mu =\bigcap \left\{ A\subset X:\overline{A}=A,\,\,\mu \in I(A) \right\}.
\end{gather*}

Consider functions of the type $\lambda:X\rightarrow [-\infty,\, 0]$. On a given set $X$ we determine a \textit{$max$-$plus$-characteristic function} $\chi_A^\oplus: X\rightarrow \mathbb{R}_{\max}$ of a subset $A\subset X$ by the rule
\begin{gather*}
^\oplus\chi_A(x)=\begin{cases}
0\, & \text{at  $x\in A$}, \\
-\infty & \text{at $x\in X\setminus A$}.
\end{cases}
\end{gather*}
For a singleton $\{x\}$ we will write $\chi_x$ instead of $\chi_{\{x\}}$.

Let $F_1$, $F_2$,\,\dots,\,$F_n$ be disjoint system of closed sets of a space $X$, and $a_1$, $a_2$,\,\dots,\,$a_n$ be non-positive real numbers. A function
\begin{gather*}
^\oplus\chi_{F_1,\,\dots,\, F_n}^{a_1,\,\dots,\, a_n}(x)=\begin{cases}
a_1\, & \text{at  $x\in F_1$}, \\
\dots,\\
a_n\, & \text{at  $x\in F_n$}, \\
-\infty & \text{at $x\in X\setminus \bigcup\limits_{i=1}^nF_n$}
\end{cases}
\end{gather*}
we call a \textit{max-plus-step-function} defined by the sets $F_1$, $F_2$,\,\dots,\, $F_n$ and the numbers $a_1$, $a_2$, \,\dots,\,$a_n$.

Note that
\begin{gather*}
^\oplus\chi_{A}^{a}(x) = a\odot\, {^\oplus\chi_{A}}(x) = \begin{cases}
0\odot a\, & \text{at  $x\in A$}, \\
-\infty & \text{at $x\in X\setminus A$}
\end{cases} = \begin{cases}
a\, & \text{at  $x\in A$}, \\
-\infty & \text{at $x\in X\setminus A$}
\end{cases}
\end{gather*}
for a set $A$ in $X$ and a non-positive number $a$. Consequently, for a disjoint system of closed sets $F_1$, $F_2$,\,\dots,\,$F_n$ in a space $X$, and non-positive real numbers $a_1$, $a_2$,\,\dots,\,$a_n$ we have
\begin{gather*}
^\oplus\chi_{F_1,\,\dots,\, F_n}^{a_1,\,\dots,\, a_n}(x) = {^\oplus\chi_{F_1}^{a_1}}(x) \oplus {^\oplus\chi_{F_2}^{a_2}}(x) \oplus\, \dots\, \oplus {^\oplus\chi_{F_n}^{a_n}}(x).
\end{gather*}

The notion of density for an idempotent measure was introduced in \cite{Akian1999}. Let $\mu\in I(X)$. Then we can define a
function $d_{\mu}\colon X\to [-\infty,\, 0]$ by the formula
\begin{equation}\label{density}
d_{\mu}(x) = \inf\{\mu(\varphi):\, \varphi\in C(X) \mbox{ such that } \varphi\le 0 \mbox{ and } \varphi(x) = 0\}, \qquad x\in X.
\end{equation}
The function $d_{\mu}$ is upper semicontinuous and is called the \textit{density} of $\mu$. Conversely, each upper semicontinuous function $f\colon X\to [-\infty, 0]$ with $\max\{f(x):\, x\in X\} = 0$ determines an idempotent measure $\nu_f$ by the formula $\nu_f(\varphi) = \bigoplus\limits_{x\in X}f(x)\odot \varphi(x)$, for $\varphi \in C(X)$.

Note that a function $f\colon X\to \mathbb{R}$ is said to be \textit{upper semicontinuous} if for each $x\in X$, and for every real number $r$ which satisfies $f(x) < r$, there exists an open neighbourhood $U\subset X$ of $x$ such that $f(x') < r$ for all $x'\in U$.

Put
\begin{multline}\label{lambda}
U_S(X)=\{\lambda\colon X\to [-\infty,\, 0]|\quad \lambda \mbox{ is upper semicontinuous and there exists a }\\
x_0\in X \mbox{ such that } \lambda(x_0)=0\}.
\end{multline}

Then we have
\begin{equation}\label{Func}
I(X) = \left\{\bigoplus\limits_{x\in X}\lambda(x)\odot \delta_x:\, \lambda\in U_S(X)\right\}.
\end{equation}

Obviously that $\bigoplus\limits_{x\in X}{^\oplus\chi_{x_0}}(x)\odot\delta_x= \delta_{x_0}$, i.~e. for a point $x_0$ by the rule (\ref{Func}) the $\max$-$\mbox{plus}$-characteristic function ${^\oplus\chi_{x_0}}$ defines the Dirac measure $\delta_{x_0}$ supported on the singleton $\{x_0\}$.

Let $A$ be a closed subset of a compactum $X$. It is easy to check that $\nu \in I(A)$ iff $\{x\in X:\, d_{\nu}(x) > -\infty\} \subset A$. Hence,
\begin{equation}\label{densupp}
\mbox{supp}\nu = \{x\in X:\, d_{\nu}(x) > -\infty\}.
\end{equation}

Consider an idempotent probability measure $\mu=\bigoplus\limits_{x\in X}\lambda(x)\odot\delta_x\in I(X)$ and a finite system $\{U_1,\, \dots,\, U_n\}$ of open sets $U_i$ such that $\mbox{supp}\mu\cap U_i\neq\varnothing$, $i=1,\, \dots,\, n$, and $\mbox{supp}\mu\subset\bigcup\limits_{i=1}^n U_i$.
Define a set
\begin{multline}\label{opensets}
\left\langle\mu;\, U_{1}\,\dots,\,U_{n};\, \varepsilon\right\rangle=\left\{\nu=\bigoplus\limits_{x\in X} \gamma(x)\odot\delta_{x}\in I(X):\, \mbox{supp}\nu\cap U_i\neq\varnothing,\,\mbox{supp}\nu\subset\bigcup\limits_{i=1}^n U_i,\right.\\
\left. \mbox{ and } |\lambda(x)-\gamma(y)|<\varepsilon \mbox{ at the points } x\in \mbox{supp}\mu \cap U_i \mbox{ and } y\in \mbox{supp}\nu \cap U_i,\, i=1,\, \dots,\, n,\right\}.
\end{multline}

\begin{Prop}
The sets of the view (\ref{opensets}) form a base of pointwise convergence topology in $I(X)$.
\end{Prop}

\begin{Proof} Let $\langle\mu;\, \varphi;\, \varepsilon\rangle$ be a prebase element, where $\varphi\in C(X)$, $\varepsilon>0$ and $\mu=\bigoplus\limits_{x\in X}\lambda(x)\odot\delta_{x}\in I(X)$. As $\varphi$ is continuous, for each point $x\in \mbox{supp}\mu$ there is its open neighbourhood $U_x$ in $X$ such that for any point $y\in U_x$ the inequality $|\varphi(x)-\varphi(y)|<\frac{\varepsilon}{2}$ holds. From the open cover $\{U_x:\, x\in \mbox{supp}\mu\}$ in $X$ of $\mbox{supp}\mu$ by owing to compactness of $\mbox{supp}\mu$ one can choose a finite subcover $\{U_i:\, i=1,\, \dots,\, n\}$. Further, for every $\nu=\bigoplus\limits_{x\in X}\gamma(x)\odot\delta_x\in \langle\mu;\, U_{1},\, \dots,\, U_{n};\, \frac{\varepsilon}{2} \rangle$ we have $|\lambda(x)-\gamma(y)| < \frac{\varepsilon}{2}$ at $x \in \mbox{supp}\mu \cap U_i$ and $y\in \mbox{supp}\nu \cap U_i$. Let us estimate the following absolute value
$|\mu(\varphi)-\nu(\varphi)|= \left|\bigoplus\limits_{x\in X}\lambda(x)\odot\varphi(x) - \bigoplus\limits_{x\in X} \gamma(x) \odot\varphi(x)\right|=a$.

Two cases are possible:

\textit{Case 1}: $\bigoplus\limits_{x\in X}\lambda(x)\odot\varphi(x) \geq \bigoplus\limits_{x\in X}\gamma_(x)\odot\varphi(x)$.  Let $\bigoplus\limits_{x\in X}\lambda(x)\odot\varphi(x)=\lambda(x')\odot\varphi(x')$. Then $x'\in U_i$ for some $i$, and
\begin{gather*}
a=\bigoplus\limits_{x\in X}\lambda(x)\odot\varphi(x) - \bigoplus\limits_{x\in X}\gamma(x)\odot\varphi(x) = \lambda(x')\odot\varphi(x') - \bigoplus\limits_{x\in X}\gamma(x)\odot\varphi(x)\le\\
\le(\mbox{for every } y\in \mbox{supp}\nu \cap U_i)\le\\
\le\lambda(x')\odot\varphi(x') - \gamma(y)\odot\varphi(y) = |\lambda(x')\odot\varphi(x') - \gamma(y)\odot\varphi(y)| \le\\
\le |\lambda(x') - \gamma(y)| + |\varphi(x') - \varphi(y)|<\frac{\varepsilon}{2} + \frac{\varepsilon}{2} = \varepsilon.
\end{gather*}

\textit{Case 2}: $\bigoplus\limits_{x\in X}\lambda(x)\odot\varphi(x) \le \bigoplus\limits_{x\in X}\gamma(x)\odot\varphi(x)$. Let $\bigoplus\limits_{x\in X}\gamma(x)\odot\varphi(x)=\gamma(x')\odot\varphi(x')$. Then $x'\in U_i$ for some $i$, and
\begin{gather*}
a=\bigoplus\limits_{x\in X}\gamma(x)\odot\varphi(x) - \bigoplus\limits_{x\in X}\lambda(x)\odot\varphi(x) = \gamma(x')\odot\varphi(x') - \bigoplus\limits_{x\in X}\lambda(x)\odot\varphi(x)\le\\
\le(\mbox{for every } y\in \mbox{supp}\mu \cap U_i)\le\\
\le\gamma(x')\odot\varphi(x') - \lambda(y)\odot\varphi(y) = |\gamma(x')\odot\varphi(x') - \lambda(y)\odot\varphi(y)| \le\\
\le |\lambda(x') - \gamma(y)| + |\varphi(x') - \varphi(y)|<\frac{\varepsilon}{2} + \frac{\varepsilon}{2} = \varepsilon.
\end{gather*}

So, in the above two cases we have $a<\varepsilon$, i.~e. $|\mu(\varphi) - \nu(\varphi)|<\varepsilon$. From here $\nu\in\langle\mu;\varphi;\varepsilon\rangle$, in other words,
\begin{gather*}
\left\langle\mu;\, U_{1},\, \dots,\, U_{n};\, \frac{\varepsilon}{2}\right\rangle \subset \langle\mu;\, \varphi;\, \varepsilon\rangle.
\end{gather*}
\end{Proof}

We recall some concepts from \cite{Pelcz1968}, and modify them for the $\max$-$\mbox{plus}$ case if necessary. Let $X$ and $Y$ be compact spaces, $f\colon X\to Y$ be a map, $f^\circ\colon C(Y)\to C(X)$ be the \textit{induced operator} defined by equality $f^\circ(\varphi) = \varphi \circ f$, $\varphi\in C(Y)$. We say that an operator $u\colon C(X)\to C(Y)$ is a $\max$-$\mbox{plus}$-\textit{linear operator} provided  $u(\alpha\odot\varphi \oplus \beta\odot\psi) = \alpha\odot u(\varphi)\oplus \beta\odot u(\psi)$ for every pair of functions $\varphi,\, \psi\in C(X)$, where $-\infty \le \alpha,\, \beta\le 0$,  $\alpha\oplus \beta = 0$. A $\max$-$\mbox{plus}$-linear operator $u\colon C(X)\to C(Y)$ is $\max$-$\mbox{plus}$-\textit{regular}  provided $\|u\| = \sup\{\|u(\varphi)\|:\, \varphi\in C(X),\, \|\varphi\|\le 1\} = 1$ and $u(1_X) = 1_Y$. A $\max$-$\mbox{plus}$-linear operator $u\colon C(X)\to C(Y)$ is said to be a $\max$-$\mbox{plus}$-\textit{linear exave} for $f$ provided $f^\circ\circ u$ is the identity on $f^\circ(C(Y))$ or equivalently $f^\circ\circ u\circ f^\circ = f^\circ$.  A $\max$-$\mbox{plus}$-\textit{regular exave} is a $\max$-$\mbox{plus}$-linear exave which is a regular operator. If $f$ is a homeomorphic embedding, then a $\max$-$\mbox{plus}$-linear exave ($\max$-$\mbox{plus}$-regular exave) for $f$ is called $\max$-$\mbox{plus}$-\textit{linear extension operator ($\max$-$\mbox{plus}$-regular extension operator)}. If $f$ is a surjective map, then a $\max$-$\mbox{plus}$-linear exave ($\max$-$\mbox{plus}$-regular exave) for $f$ is called $\max$-$\mbox{plus}$-\textit{linear averaging operator ($\max$-$\mbox{plus}$-regular averaging operator)}.

Remind, in category theory a \textit{monomorphism} (an \textit{epimorphism}) is a left-cancellative (respectively, right-cancellative) morphism, that is, a morphism $f \colon Z \to X$  (respectively, $f \colon X \to Y$) such that, for each pair of  morphisms  $g_1$, $g_2\colon Y \to Z$ the following implication holds
\begin{gather*}
f \circ g_1 = f \circ g_2 \Rightarrow g_{1}=g_{2} \qquad (\mbox{respectively, }\, g_{1}\circ f=g_{2}\circ f\Rightarrow g_{1}=g_{2}).
\end{gather*}

If $u$ is an exave for $f\colon X\to Y$ and $y\in f(X)$, then for every function $\varphi\in C(Y)$ we have
\begin{equation}
(u\circ f^{\circ})(\varphi)(y) = \varphi(y).
\end{equation}

\begin{Prop}\label{ucircf}
Let $f\colon X\to Y$ be a map. A $\max$-$\mbox{\rm plus}$-regular operator $u\colon C(X)\to C(Y)$ is a $\max$-$\mbox{\rm plus}$-regular extension (respectively, averaging) operator if and only if $f^{\circ}\circ u = \mbox{\rm id}_{C(X)}$ (respectively, $u\circ f^{\circ} = \mbox{\rm id}_{C(Y)}$).
\end{Prop}
\begin{Proof}
Let $u$ be a $\max$-$\mbox{\rm plus}$-regular extension (respectively, averaging) operator. Then $f^{\circ}\colon C(Y)\to C(X)$ is an epimorphism (respectively, monomorphism). Thence $f^{\circ} \circ u\circ f^{\circ} = f^{\circ} =\mbox{id}_{C(X)} \circ f^{\circ}$ implies $f^{\circ}\circ u = \mbox{id}_{C(X)}$ (respectively, $f^{\circ} \circ u\circ f^{\circ} = f^{\circ} = f^{\circ}\circ \mbox{id}_{C(Y)}$ implies $u \circ f^{\circ}= \mbox{id}_{C(Y)}$).

Let $u$ be a $\max$-$\mbox{\rm plus}$-regular operator and $f^{\circ}\circ u = \mbox{id}_{C(X)}$. It requires to show $f\colon X\to Y$ is an embedding. Suppose $f(x_1) = f(x_2)$, $x_1$, $x_2\in X$. Assume there exists a function $\varphi\in C(X)$ such that $\varphi(x_1) \neq \varphi(x_2)$. Conversely, we have $\varphi(x_1) = f^{\circ}\circ u (\varphi)(x_1) = u (\varphi)(f(x_1)) = u (\varphi)(f(x_2)) = f^{\circ}\circ u (\varphi)(x_2) = \varphi(x_2)$. We get a contradiction. So, $x_1 = x_2$.

Let $u$ be a $\max$-$\mbox{\rm plus}$-regular operator and $u \circ f^{\circ} = \mbox{id}_{C(Y)}$. We should show that $f\colon X\to Y$ is a surjective map. Suppose $f$ is not so. Then $Y\setminus f(X) \ne \varnothing$ and for every $y\in Y\setminus f(X)$, since the image $f(X)$ is a compact space, any $\varphi\colon f(X)\to \mathbb{R}$ has different extensions $\varphi_1$, $\varphi_2\colon Y\to \mathbb{R}$ such $\varphi_1(y) \ne \varphi_2(y)$. Hence, $\varphi_1\ne \varphi_2$. On the other hand $\varphi_1 = u \circ f^{\circ} (\varphi_1) = u \circ f^{\circ} (\varphi_2) = \varphi_2$. The got contradiction finishes the proof.

\end{Proof}

An epimorphism $f\colon X\to Y$ is said to be a $\max$-$\mbox{plus}$-\textit{Milutin epimorphism}  provided it permits a
$\max$-$\mbox{plus}$-regular averaging operator. A compact space $X$ is a $\max$-$\mbox{plus}$-\textit{Milutin space}  if there exists a $\max$-$\mbox{plus}$-Milutin epimorphism $f\colon D^{\tau}\to X$ \cite{Pelcz1968}. Every compactum is a Milutin space (\cite{FedorchukGenTop}, Corollary VIII.4.6.). Analogously, every compactum is a $\max$-$\mbox{plus}$-Milutin space.\\

\section{An analog of the Uspenskii's metrics}

Every zero-dimensional space of the weight $\mathfrak{m}\ge \aleph_0$ embeds into Cantor cube $D^{\mathfrak{m}}$.  Consequently, a zero-dimensional compactum is a $\max$-$\mbox{plus}$-Milutin space.

Let $\mu_1 = \bigoplus\limits_{x\in X}\lambda_1(x)\odot\delta_x$, $\mu_2 = \bigoplus\limits_{x\in X}\lambda_2(x)\odot\delta_x\, \in I(X)$.  Put
\begin{gather*}
\Lambda_{1\,2} = \Lambda (\mu_1,\, \mu_2) = \{\xi \in I(X^2):\, I(\pi_i)(\xi) = \mu_i,\, i= 1,\, 2\},
\end{gather*}
where $\pi_i \colon\, X \times  X \to X$ is the projection onto $i$-th factor, $i = 1,\, 2$. We will show the set
$\Lambda (\mu_1,\, \mu_2)$ is nonempty. Let $x_{i\, 0}\in \mbox{supp}\mu_i$ be points such that $\lambda_i(x_{i\,0}) = 0$, $i=1,\,2$. Then the directly checking shows that $I(\pi_i)(\xi) = \mu_i$, $i = 1,\, 2$, for all $\xi \in I(X^2)$ of the form $\xi  = \xi^0 \oplus R(\mu_1, \mu_2)$. Here
\begin{gather*}
\xi^0 = 0\odot \delta_{(x_{1\,0},\, x_{2\,0})} \bigoplus\limits_{x\in X\setminus \{x_{1\,0}\}} \lambda_2(x)\odot\delta_{(x_{1\,0},\, x)}\oplus \bigoplus\limits_{x\in X\setminus \{x_{2\,0}\}} \lambda_1(x)\odot\delta_{(x,\, x_{2\,0})}
\end{gather*}
is an idempotent probability measure on $X^2$, and
\begin{gather*}
R(\mu_1,\, \mu_2) = \bigoplus\limits_{\substack{x\in X\setminus \{x_{1\,0}\}\\ y\in X\setminus \{x_{2\,0}\}}} \gamma(x,\, y) \odot \delta_{(x,\, y)}
\end{gather*}
is some functional on $C(X)$ where
\begin{gather*}
-\infty\le \gamma(x,\, y)\le \min\{\lambda_1(x),\, \lambda_2(y)\}, \, x\in M,\, y\in N,\, M\subset X\setminus \{x_{1\,0}\},\, N\subset X\setminus \{x_{2\,0}\}.
\end{gather*}

Thus $\xi \in  \Lambda(\mu_1,\, \mu_2)$, i.~e. $\Lambda(\mu_1,\, \mu_2) \neq \varnothing$. In fact, here more is proved: it is easy to see if $|X|\ge 2$ and $|Y|\ge 2$
then quantity of the numbers $\gamma(x,\, y)$ is uncountable. From here one concludes that the potency of the set $\Lambda(\mu_1,\, \mu_2)$ is no less than continuum potency as soon as each of the supports $\mbox{supp}\mu_i$, $i=1,\,2$, contains no less than two points.

Note that $\xi = \xi^0$ if one takes empty set as $K$ and $M$.

Idempotent probability measures $\xi\in I(X^2)$ with $I(\pi_i)(\xi) = \mu_i,\, i= 1,\, 2$ we will call as $(\mu_1,\,\mu_2)$-\textit{admissible  measures}.

The following statement is rather evident.
\begin{Prop}\label{prfunct}
Let $\mu_i = \bigoplus\limits_{x\in X} \lambda_i(x)\odot \delta_x$, $i=1,\,2$, be idempotent probability measures. Then every $(\mu_1,\,\mu_2)$-admissible measure $\xi = \bigoplus\limits_{(x,\, y)\in X^2} \lambda_{1\,2}(x,\, y)\odot \delta_{(x,\, y)} \in I(X^2)$ satisfies the following equalities:
\begin{gather*}
\lambda_1(x) =  \bigoplus\limits_{y\in X} \lambda_{1\,2}(x,\, y),\quad x\in X,\quad \mbox{ and }\quad \lambda_2(y) =  \bigoplus\limits_{x\in X} \lambda_{1\,2}(x,\, y),\quad y\in X.
\end{gather*}
\end{Prop}

Consider a compactum $(X,\, \rho)$. We define a function $d_{I}\colon I(X)\times I(X)\to \mathbb{R}$ by the formula
\begin{gather*}
d_{I}(\mu_1,\,\mu_2) = \inf\{\xi(\rho):\, \xi\in \Lambda_{1\,2}\}.
\end{gather*}

This function was offered by V.~V.~Uspenskii and in \cite{Fedorchuk1990} it was proved that it is a metrics on the space of probability measures.  Its analog for idempotent probability measures is not metrics on the space of idempotent probability measures.

\begin{Prop}\label{admexsist}
For every pair $\mu_1,\, \mu_2\in I(X)$ there exists a $(\mu_1,\,\mu_2)$-admissible idempotent probability measure $\xi\in I(X^2)$ such that
\begin{gather*}
d_{I}(\mu_1,\,\mu_2) = \xi(\rho).
\end{gather*}
\end{Prop}
\begin{Proof}
Consider a sequence $\{\xi_n\}$ of $(\mu_1,\,\mu_2)$-admissible idempotent probability measures such that $\xi_n(\rho)\longrightarrow d_{I}(\mu_1,\,\mu_2)$. Passing in case of need to a subsequence, owing to compactness of $I(X^2)$, it is possible to assume that $\{\xi_n\}$ tends to some $\xi\in I(X^2)$. Since the projections $I(\pi_i)$ are continuous, $\xi$ is $(\mu_1,\,\mu_2)$-admissible.
Further, for an arbitrary $\varepsilon>0$ there exists $n_0$ such that $\xi_n\in \langle\xi;\, \rho;\, \varepsilon\rangle$ for all $n\ge n_0$, where $\langle\xi;\, \rho;\, \varepsilon\rangle$ is a prebase neighbourhood of $\xi$ in the pointwise convergence topology on $I(X^2)$. So, $|\xi(\rho)-\xi_n(\rho)|<\varepsilon$. Consequently, $d_{I}(\mu_1,\,\mu_2) = \xi(\rho)$.

\end{Proof}

\begin{Prop}\label{pseudomet}
The function $d_{I}$ is a pseudometric on $I(X)$.
\end{Prop}
\begin{Proof} Since each $\xi\in I(X^2)$ is order-preserving then the inequality  $\rho\ge 0$ immediately implies $d_{I}\ge 0$. So, $d_{I}$ is nonnegative. Obviously, $d_{I}$ is symmetric.

Let $\mu_1= \mu_2 = \mu$. There exists $\lambda\in U_S(X)$ such that $\mu = \bigoplus\limits_{x\in X}\lambda(x)\odot\delta_x$. Then $\xi_{\mu} = \bigoplus\limits_{x\in X}\lambda(x)\odot\delta_{(x,\, x)}$ is a  $(\mu_1,\,\mu_2)$-admissible idempotent probability measure, and
\begin{gather*}
0\le d_{I}(\mu_1,\,\mu_2) = \inf\{\xi(\rho):\, \xi\in \Lambda_{1\,2}\}\le \xi_{\mu} (\rho) = \bigoplus\limits_{x\in X}\lambda(x) = 0,
\end{gather*}
i.~e. $d_{I}(\mu^1,\,\mu^2) = 0$.

Let us show that the triangle inequality is true as well. Take arbitrary triple $\mu_i\in I(X)$, $i=1,\,2,\,3$. Let $\mu_{1\,2},\,\mu_{2\,3}\in I(X^2)$ be $(\mu_1,\, \mu_2)$- and $(\mu_2,\, \mu_3)$-admissible measures such that $d_{I}(\mu_1,\, \mu_2) = \mu_{1\,2}(\rho)$ and $d_{I}(\mu_2,\, \mu_3) = \mu_{2\,3}(\rho)$, respectively. Put
\begin{gather*}
X_1 = X_2 = X_3 = X,\qquad X_{1\,2\,3} = X^3 = X_1\times X_2\times X_3,\qquad X_{i\,j} = X^2 = X_i\times X_j,
\end{gather*}
and let
\begin{gather*}
\pi^{1\,2\,3}_{i\,j}\colon X_{1\,2\,3}\to X_{i\,j},\qquad  \pi^{i\,j}_{k}\colon X_{i\,j}\to X_k,\qquad 1\le i< j\le 3,\qquad k\in \{i,\,j\},
\end{gather*}
be corresponding projection.

According to Corollary 4.3 \cite{Zar2010} the functor $I$ is bicommutative. Using this fact one can similarly to Lemma 4 \cite{Fedorchuk1990} show that for idempotent probability measures
\begin{gather*}
\mu_2 \in I(X_2), \qquad \mu_{1\,2}\in I(X_{1\,2}),\qquad \mu_{2\,3}\in I(X_{2\,3})
\end{gather*}
such that
\begin{gather*}
I(\pi^{1\,2}_{2})(\mu_{1\,2}) = \mu_2 = I(\pi^{2\,3}_{2})(\mu_{2\,3}),
\end{gather*}
there exists $\mu_{1\,2\,3}\in I(X_{1\,2\,3})$ which satisfies the equalities
\begin{gather*}
I(\pi^{1\,2\,3}_{1\,2})(\mu_{1\,2\,3}) = \mu_{1\,2}\qquad \mbox{ and }\qquad I(\pi^{1\,2\,3}_{2\,3})(\mu_{1\,2\,3}) = \mu_{2\,3}.
\end{gather*}

Set $\mu_{1\,3} = I(\pi^{1\,2\,3}_{1\,3})(\mu_{1\,2\,3})$. Then according to Proposition \ref{prfunct}  $\mu_{1\,3}$ is a $(\mu_1,\, \mu_3)$-admissible idempotent probability measure. Using Proposition \ref{prfunct}, we obtain
\begin{gather*}
d_{I}(\mu_1,\, \mu_2) +d_{I}(\mu_2,\, \mu_3) =  \mu_{1\,2}(\rho) + \mu_{2\,3}(\rho) =\\
=\bigoplus_{(x_1,\,x_2)\in X_{1\,2}} d_{\mu_{1\,2}}(x_1,\,x_2)\odot \rho (x_1,\,x_2) + \bigoplus_{(x_2,\,x_3)\in X_{2\,3}} d_{\mu_{2\,3}}(x_2,\,x_3)\odot \rho (x_2,\,x_3) =\\
=\bigoplus_{(x_1,\,x_2,\,x_3)\in X_{1\,2\,3}} d_{\mu_{1\,2\,3}}(x_1,\,x_2,\,x_3)\odot \rho (x_1,\,x_2) + \bigoplus_{(x_1,\,x_2,\,x_3)\in X_{1\,2\,3}} d_{\mu_{1\,2\,3}}(x_1,\,x_2,\,x_3)\odot \rho (x_2,\,x_3) \ge\\
\ge \bigoplus_{(x_1,\,x_2,\,x_3)\in X_{1\,2\,3}} \left(d_{\mu_{1\,2\,3}}(x_1,\,x_2,\,x_3)\odot \rho (x_1,\,x_2) + d_{\mu_{1\,2\,3}}(x_1,\,x_2,\,x_3)\odot \rho (x_2,\,x_3)\right) =\\
= \bigoplus_{(x_1,\,x_2,\,x_3)\in X_{1\,2\,3}}d_{\mu_{1\,2\,3}}(x_1,\,x_2,\,x_3)\odot \left(\rho (x_1,\,x_2) + \rho (x_2,\,x_3)\right)\ge\\
\ge \bigoplus_{(x_1,\,x_2,\,x_3)\in X_{1\,2\,3}}d_{\mu_{1\,2\,3}}(x_1,\,x_2,\,x_3)\odot \rho (x_1,\,x_3) = \\
=\bigoplus_{(x_1,\,x_3)\in X_{1\,3}}d_{\mu_{1\,3}}(x_1,\,x_3)\odot \rho (x_1,\,x_3) = \mu_{1\,3}(\rho) \ge d_{I}(\mu_1,\,\mu_3),
\end{gather*}
i.~e. $d_{I}(\mu_1,\,\mu_3)\le d_{I}(\mu_1,\, \mu_2) +d_{I}(\mu_2,\, \mu_3)$. Here $d_\nu$ is the density function of the corresponding measure $\nu$ (see page \pageref{density}).

\end{Proof}

Unlike usual probability measures, the function $d_{I}$ is not a metrics.

\begin{Ex}\label{exampseudometr}
{\rm Let $(X,\, \rho)$ be a metric space, $x,\,y\in X$ be points such that $\rho(x,\, y)=1$.  Consider idempotent probability measures $\mu_1=0\odot\delta_x\oplus (-2)\odot\delta_y$ and $\mu_2=0\odot\delta_x\oplus (-4)\odot\delta_y$. One can directly check that the idempotent probability measure $\xi=0\odot \delta_{(x,\,x)}\oplus (-2)\odot \delta_{(y,\,x)}\oplus (-4)\odot \delta_{(x,\,y)}$ is $(\mu_1,\, \mu_2)$-admissible, and $\xi(\rho)=0$. That is why $d_{I}(\mu_1,\, \mu_2)=0$, though $\mu_1\neq \mu_2$.}
\end{Ex}

Example \ref{exampseudometr} shows that the functors $P$ of probability measures and $I$ of idempotent probability measures are not isomorphic.\\

\section{On a metrics on the space of idempotent probability measures}

Let $(X,\, \rho)$ be a metric compact space. We suggest a distance function $\rho_I\colon I(X)\times I(X) \to \mathbb{R}$ as follows
\begin{equation}\label{rhoI}
\rho_I(\mu_1,\, \mu_2) = \inf\{\sup\{\xi(\rho) \oplus \rho(x,\, y):\, (x,\, y)\in \mbox{supp}\xi\}:\, \xi\in \Lambda_{1\,2}\}.
\end{equation}

\begin{Th}\label{met}
The function $\rho_I$ is a metrics on $I(X)$  which is an extension of the metric $\rho$.
\end{Th}
\begin{Proof}
Obviously, $\rho_I$ is nonnegative and symmetric. If $\mu_1=\mu_2$ then similarly to the proof of Proposition \ref{pseudomet} one can show that $\rho_I (\mu_1, \mu_2) = 0$. Inversely, let $\rho_I (\mu_1, \mu_2) = 0$. Then it there exist a $\xi\in \Lambda_{1\, 2}$ such that $\rho(x,\, y) = 0$ for all $(x,\, y)\in \mbox{supp}\xi$. Consequently $\mbox{supp}\xi$ must lie in the diagonal $\Delta(X) = \{(x,\, x):\, x\in X\}$. Applying Proposition \ref{prfunct}, we have $d_{\mu_1} = d_{\mu_2}$, which implies $\mu_1 = \mu_2$. It remains to check the triangle axiom. But the checking consists only of the repeating of procedure at the proof of Proposition \ref{pseudomet}.

For every pair of Dirac measures $\delta_x$, $\delta_y$,  $x,\, y\in X$, the uniqueness of $(\delta_x,\, \delta_y)$-admissible measure $\xi\in I(X^2)$, $\xi = 0\odot \delta_{(x,\,y)}$, implies that
\begin{gather*}
\rho_{I}(\delta_x,\, \delta_y) = \xi(\rho) \oplus \rho(x,\, y) = 0\odot \delta_{(x,\,y)}(\rho) \oplus \rho(x,\, y) = \rho(x,\,y).
\end{gather*}
From here we get that $\rho_{I}$ is an extension of $\rho$.

\end{Proof}

Since every idempotent probability measure is order-preserving, from the construction of the metrics  $\rho_{I}$ we obtain the following statement.

\begin{Prop}{\rm
$\mbox{diam}(I(X),\, \rho_{I}) = \mbox{diam}(X,\, \rho)$.}
\end{Prop}
\begin{Proof}
Indeed, we have $0\le \xi(\rho) \le \mbox{diam}(X,\, \rho)$ as $0\le \rho \le \mbox{diam}(X,\, \rho)$.

\end{Proof}

\begin{Prop}\label{minusinfp}
Let $X$ be a compactum and a sequence $\{\mu_n\}\subset I(X)$ converge to $\mu_0\in I(X)$ with respect to pointwise convergence topology. Then for every open neighbourhood $U$ of the diagonal $\Delta(X)=\{(x,\, x):\, x\in X\}$ there exist a positive integer $n$ and a $(\mu_0,\, \mu_n)$-admissible measure $\mu_{0\, n}\in I(X^2)$ such that
\begin{equation}\label{minusinfe}
\bigoplus\limits_{(x,\, y)\in X^2\setminus U} d_{\mu_{0\, n}}(x,\, y)\odot\rho(x,\, y) = -\infty.
\end{equation}
\end{Prop}
\begin{Proof}
At first we consider the case of zero-dimensional compactum $X$. There exists a disjoint clopen cover $\{V_1,\, \dots,\, V_n\}$ of $X$ (i.~e. a cover, which consists of open-closed sets of $X$) such that $V_i\times V_i\subset U$ for each $i = 1,\, \dots,\, n$. As $\mu_n\rightarrow \mu$ there exists $n$ such that $\mu_n\in \langle\mu;\, ^\oplus\chi_{V_1},\, ^\oplus\chi_{V_2},\, \dots,\, ^\oplus\chi_{V_n};\, \varepsilon\rangle$. We will determine $(\mu_0,\, \mu_n)$-admissible measure $\mu_{0\, n}\in I(X^2)$.

There exists a base of the compactum $X$ consisting of clopen sets
\begin{gather*}
V^{\varepsilon_1\varepsilon_2\dots\varepsilon_k}_i, \qquad 1\le i\le s,\qquad \varepsilon_k\in\{0,\, 1\},\qquad 1\le k <\infty,
\end{gather*}
such that
\begin{itemize}
\item[$1)$] $V^0_i\cup V^1_i = V_i$;
\item[$2)$] $V^0_i\cap V^1_i = \varnothing$;
\item[$3)$] $V^{\varepsilon_1\varepsilon_2\dots\varepsilon_k 0}_i \cup V^{\varepsilon_1\varepsilon_2\dots\varepsilon_k 1}_i = V^{\varepsilon_1\varepsilon_2\dots\varepsilon_k}_i$;
\item[$4)$] $V^{\varepsilon_1\varepsilon_2\dots\varepsilon_k 0}_i \cap V^{\varepsilon_1\varepsilon_2\dots\varepsilon_k 1}_i = \varnothing$.
\end{itemize}

The sets $V^{\varepsilon_1\varepsilon_2\dots\varepsilon_k}_i\times V^{\varepsilon'_1\varepsilon'_2\dots\varepsilon'_k}_{i'}$ form a base of the compactum $X_{1\, 2}$. To determine $\mu_{0\, n}$ it is enough to construct its density function. Let $\mu_0=\bigoplus\limits_{x\in X} \lambda_0(x)\odot \delta_x$, $\mu_n=\bigoplus\limits_{x\in X} \lambda_n(x)\odot \delta_x$. We set
\begin{gather*}
\lambda^{\varepsilon_1\dots\varepsilon_k,\,\varepsilon'_1\dots\varepsilon'_k}_{i\,i'} = \bigoplus\limits_{(x,\, y)\in X\times X}(\lambda_0(x)\odot \lambda_n(y))\odot \delta_{(x,\, y)}(^\oplus\chi_{V^{\varepsilon_1\dots\varepsilon_k}_i\times V^{\varepsilon'_1\dots\varepsilon'_k}_{i'}}),
\end{gather*}
i.~e.
\begin{gather*}
\lambda^{\varepsilon_1\dots\varepsilon_k,\,\varepsilon'_1\dots\varepsilon'_k}_{i\,i'} = \bigoplus\limits_{(x,\, y)\in V^{\varepsilon_1\dots\varepsilon_k}_i\times V^{\varepsilon'_1\dots\varepsilon'_k}_{i'}}\lambda_0(x)\odot \lambda_n(y).
\end{gather*}

It is clear that
\begin{gather*}
\lambda^{\varepsilon'_1\dots\varepsilon'_k}_{i'} = \bigoplus\limits_{i = 1}^{s} \lambda^{\varepsilon_1\dots\varepsilon_k,\,\varepsilon'_1\dots\varepsilon'_k}_{i\,i'}\qquad \mbox{and} \qquad
\lambda^{\varepsilon_1\dots\varepsilon_k}_{i} = \bigoplus\limits_{i' = 1}^{s} \lambda^{\varepsilon_1\dots\varepsilon_k,\,\varepsilon'_1\dots\varepsilon'_k}_{i\,i'},
\end{gather*}
where
\begin{gather*}
\lambda^{\varepsilon_1\dots\varepsilon_k}_{i} = \bigoplus\limits_{x\in X}\lambda_0(x)\odot \delta_x(^\oplus\chi_{V^{\varepsilon_1\dots\varepsilon_k}_i}) = \bigoplus\limits_{x\in V^{\varepsilon_1\dots\varepsilon_k}_i}\lambda_0(x)
\end{gather*}
and
\begin{gather*}
\lambda^{\varepsilon'_1\dots\varepsilon'_k}_{i'} = \bigoplus\limits_{x\in X}\lambda_n(x)\odot \delta_x(^\oplus\chi_{V^{\varepsilon'_1\dots\varepsilon'_k}_i}) = \bigoplus\limits_{x\in V^{\varepsilon'_1\dots\varepsilon'_k}_{i'}}\lambda_n(x).
\end{gather*}

Put
\begin{gather*}
d_{\mu_{0\, n}} = \lim\limits_{s\to \infty} \bigoplus\limits_{i,\, i'=1}^{s} {^\oplus\chi}^{\lambda^{\varepsilon_1\dots\varepsilon_k,\,\varepsilon'_1\dots\varepsilon'_k}_{i\,i'}}_ {V^{\varepsilon_1\dots\varepsilon_k}_i\times V^{\varepsilon'_1\dots\varepsilon'_k}_{i'}}.
\end{gather*}
Then $d_{\mu_{0\, n}}$ is an upper semicontinuous function on $X^2$ and  $\mu_{0,\, n} = \bigoplus\limits_{(x,\, y)\in X^2}d_{\mu_{0\, n}}(x,\, y)\odot \delta_{(x,\, y)}$ is a $(\mu_0,\, \mu_n)$-admissible measure with $\mbox{supp}\mu_{0,\, n} \subset U$. Consequently,
$\bigoplus\limits_{(x,\, y)\in X^2\setminus U} d_{\mu_{0\, n}}(x,\, y) = -\infty$ and, the equation (\ref{minusinfe}) is proved for the zero-dimensional case.

Now let $X$ be an arbitrary compactum. There exists a zero-dimensional compactum $Z$, a $\max$-$\mbox{plus}$-Milutin epimorphism $f\colon Z\to X$ and a $\max$-$\mbox{plus}$-regular averaging operator $u\colon C(Z)\to C(X)$ corresponding to this epimorphism. The dual $\max$-$\mbox{plus}$-map $u^{\oplus}$ which we define by the equality $u^{\oplus}(\mu)(\varphi) = \mu(u(\varphi))$, $\varphi\in C(Z)$, generates an embedding $u^{\oplus}\colon I(X)\to I(Z)$.

For idempotent probability measures $\mu'_0 = u^{\oplus} (\mu_0)$ and $\mu'_n = u^{\oplus} (\mu_n)$ there exists  $(\mu'_0,\, \mu'_n)$-admissible idempotent probability measure $\mu'_{0,\, n} = \bigoplus\limits_{(x',\, y')\in Z^2}d_{\mu'_{0\, n}}(x',\, y')\odot \delta_{(x',\, y')}\in I(Z\times Z)$ such that
\begin{gather*}
\bigoplus\limits_{(x',\, y')\in Z^2\setminus (f\times f)^{-1}(U)} d_{\mu'_{0\, n}}(x',\, y')\odot\rho(x',\, y') = -\infty.
\end{gather*}
Put $\mu_{0,\, n} = I(f\times f)(\mu'_{0\, n})$. Then for every $\varphi\in C(X^2)$ we have
\begin{multline*}\mu_{0,\, n}(\varphi) = I(f\times f)(\mu'_{0\, n})(\varphi) = \mu'_{0\, n} (\varphi\circ (f\times f)) = \bigoplus\limits_{(x',\, y')\in Z^2} d_{\mu'_{0\, n}}(x',\, y')\odot \varphi \circ (f\times f)(x',\, y') = \\
= \bigoplus\limits_{(x',\, y')\in Z^2} d_{\mu'_{0\, n}}(x',\, y')\odot \varphi (f(x'),\, f(y')) =
\bigoplus\limits_{(x,\, y)\in X^2} d_{\mu'_{0\, n}}(x,\, y))\odot \delta_{(x,\, y)}(\varphi),
\end{multline*}
i.~e. $\mu_{0,\, n} = \bigoplus\limits_{(x,\, y)\in X^2} d_{\mu'_{0\, n}}(x,\, y))\odot \delta_{(x,\, y)}$. Here $d_{\mu'_{0\, n}}(x,\, y) = \bigoplus\limits_{(x',\, y')\in (f\times f)^{-1}(x,\, y)} d_{\mu'_{0\, n}}(x',\, y')$. That is why
\begin{gather*}
\bigoplus\limits_{(x,\, y)\in X^2\setminus U} d_{\mu'_{0\, n}}(x,\, y)\odot\rho(x,\, y) = -\infty.
\end{gather*}
So, $\mu_{0,\, n} = I(f\times f)(\mu'_{0\, n})$ satisfies (\ref{minusinfe}). It remains to show that $\mu_{0,\, n}$ is $(\mu_0,\, \mu_n)$-admissible.

A diagram
\begin{equation}
\begin{CD}\label{comdiag}
Z\times Z @>{f\times f}>> X\times X \\
@VV{\theta^{12}_1}V @VV{\pi^{12}_1}V \\
Z @>{f}>> X
\end{CD}
\end{equation}
is commutative, where $\theta^{12}_1$, $\pi^{12}_1$ are projections onto the first corresponding factors. Then
\begin{multline*}
I(\pi^{12}_1)(\mu_{0\, n}) =  I(\pi^{12}_1)\circ I(f\times f)(\mu'_{0\, n}) = I(\pi^{12}_1\circ (f\times f))(\mu'_{0\, n}) = \\
=(\mbox{owing to commutativity of the diagram (\ref{comdiag})}) =\\
= I(f\circ\theta^{12}_1)(\mu'_{0\, n}) = I(f)\circ I(\theta^{12}_1)(\mu'_{0,\, n}) = I(f)(\mu'_0) = I(f)(u^\oplus(\mu_0)),
\end{multline*}
i.~e. for every $\varphi\in C(X)$ we have
\begin{multline*}
I(\pi^{12}_1)(\mu_{0\, n})(\varphi) = I(f)(u^\oplus(\mu_0)) (\varphi) = u^\oplus(\mu_0) (\varphi \circ f) = u^\oplus(\mu_0) (f^{\circ}(\varphi)) = \mu_0(u\circ f^{\circ}(\varphi)) = \\
= (\mbox{with respect to Proposition \ref{ucircf}}) = \mu_0(\varphi).
\end{multline*}
Thus, $I(\pi^{12}_1)(\mu_{0\, n}) = \mu_0$. Similarly, $I(\pi^{12}_2)(\mu_{0\, n}) = \mu_n$. The Proposition is proved.

\end{Proof}

\begin{Th}
The metrics $\rho_{I}$ generates pointwise convergence topology on $I(X)$.
\end{Th}
\begin{Proof}
Let $\{\mu_n\} \subset I(X)$ be a sequence and $\mu_0 \in I(X)$. Suppose the sequence converges to $\mu_0$ with respect to the pointwise convergence topology but not by $\rho_I$. Passing in case of need to a subsequence, it is possible to regard that
\begin{gather*}
\rho_I(\mu_n,\, \mu_0)\ge a> 0\qquad \mbox{ for all positive integer }\quad n.
\end{gather*}
Consider an open neighbourhood of the diagonal $\Delta(X)$:
\begin{gather*}
U = \left\{(x,\, y)\in X^2:\, \rho(x,\, y)<\frac{a}{2}\right\}.
\end{gather*}
By virtue of Proposition \ref{minusinfp} there exist a positive integer $n$ and a $(\mu_0,\, \mu_n)$-admissible measure $\mu_{0\, n}\in I(X^2)$ such that
\begin{gather*}
\bigoplus\limits_{(x,\, y)\in X^2\setminus U} d_{\mu_{0\, n}}(x,\, y)\odot\rho(x,\, y) = -\infty.
\end{gather*}
Therefore, $\mbox{supp}\mu_{0\, n}\subset U$, and
\begin{multline*}
\rho_I(\mu_n,\, \mu_0) \le \sup\left\{\mu_{0\, n} (\rho)\oplus \rho(z,\, t):\, (z,\, t)\in \mbox{supp}\mu_{0\, n} \right\} = \\
=\sup\limits_{(z,\, t)\in {\scriptsize\mbox{supp}}\mu_{0\, n}}\left\{\left(\bigoplus\limits_{(x,\, y)\in X^2}d_{\mu_{0\, n}}(x,\, y)\odot\rho(x,\, y)\right)\oplus \rho(z,\, t)\right\}=\\
=\sup\limits_{(z,\, t)\in {\scriptsize\mbox{supp}}\mu_{0\, n}}\left\{\left(\bigoplus\limits_{(x,\, y)\in X^2\setminus U}d_{\mu_{0\, n}}(x,\, y)\odot\rho(x,\, y) \oplus \sup\limits_{(x,\, y)\in U}d_{\mu_{0\, n}}(x,\, y)\odot\rho(x,\, y)\right)\oplus \rho(z,\, t)\right\} = \\
= \sup\limits_{(z,\, t)\in {\scriptsize\mbox{supp}}\mu_{0\, n}}\left\{\left(\sup\limits_{(x,\, y)\in U}d_{\mu_{0\, n}}(x,\, y)\odot\rho(x,\, y)\right)\oplus \rho(z,\, t)\right\} \le\\
\le \sup\limits_{(z,\, t)\in U}\left\{\left(\sup\limits_{(x,\, y)\in U}d_{\mu_{0\, n}}(x,\, y)\odot\rho(x,\, y)\right)\oplus \rho(z,\, t)\right\} \le\frac{a}{2} < a.
\end{multline*}
The obtained contradiction finishes the proof.

\end{Proof}

\section*{Acknowledgements}

\end{document}